\documentclass[12pt,a4paper]{article}
\usepackage[top=3.5cm, bottom=3.5cm, left=3cm, right=3cm]{geometry}
\usepackage[latin1]{inputenc}
\usepackage{csquotes}
\usepackage{amsmath}
\usepackage{amsthm}
\usepackage{amsfonts}
\usepackage{amssymb}
\usepackage{graphics}
\usepackage{float}
\usepackage[hang,flushmargin]{footmisc} 

\usepackage{amsmath,amsthm,amssymb,graphicx, multicol, array}
\usepackage{enumerate}
\usepackage{enumitem}
\usepackage{xcolor}

\usepackage[pagebackref,bookmarks,colorlinks,breaklinks]{hyperref}

\hypersetup{linkcolor=blue,citecolor=red,filecolor=blue,urlcolor=blue} 

\usepackage{tabto}

\numberwithin{equation}{section}

\usepackage{tikz}
\usepackage{pgfplots}

\usetikzlibrary{matrix}

\usepackage{mathrsfs}

\DeclareMathOperator{\ord}{ord}

\newtheorem{thm}{Theorem}[section]
\newtheorem{lem}{Lemma}[section]

\newtheorem{cor}{Corollary}[section]
\newtheorem{dfn}{Definition}[section]

\newcommand{\F}{\mathbb{F}}

\usepackage{fancyhdr}
\title{Formulas for the Square Roots Mod $p$}
\date{}
\author{N. A. Carella}
\begin{document}
	
\maketitle	
\begin{abstract}
A method of constructing specific polynomial representations $f(x) \in \F_p[x]$ of the square roots function
modulo a prime $p = 2^kn + 1$, $n$ odd, is presented. The formulas for the cases $ k = 2, 3$ and $ 4$ are given. \let\thefootnote\relax\footnote{ \today \date{} 
	}
\end{abstract}

\section{Introduction}\label{S8180}
The polynomial representations of various functions on finite fields are important components in modern information science. The discrete logarithm, discrete exponentiation, and $n$th root functions are of significant interest in the design of cryptographic protocols. For example, the computations of the square root modulo a prime of the form $p = 2^r - 2^s + 1$ is a
step in the standard elliptic curve cryptographic protocol.
A polynomial $f(x) \in \F_p[x]$ is called a polynomial representation of the square root function mod p if it satisfies the equation
$\sqrt{x}\equiv  \pm f (x)\bmod p$ whenever $x\ne0$ is a quadratic residue.
Polynomial interpolation is the main tool used to construct polynomial representations of functions on finite fields. This is
a time tested method and works in every case. Polynomial interpolation is used in \cite{AN2003} to prove the existence of polynomial
representations $f(x)$ of the square root function of degree $\deg(f) \leq (p-3)/2$ and length (the number of nonzero terms) at
most $2^k-1$. A few specific polynomial representations are also given.
In this note a different method will be utilized to generate specific polynomial representations of the square root function
modulo $p$. The result also improves the degree estimate. The first few of these polynomials are also computed.	
	
\section{Foundation}\label{S8181}

\begin{dfn} \label{dfn8181.200D} \hypertarget{dfn8181.200D}  Let $p = 2^kn + 1$ be prime, where $k\geq3$ and $n\geq1$ is odd, and let $z$ be a quadratic nonresidue modulo $p$. The group of $2^k$th root of unity in the finite field $\F_p$ is denoted by
	$$\mu_{2^k}=\{1,z^n,z^{2n},z^{3n},\ldots,z^{(2^{k}-1)n}\}.$$
The group of roots of unity is a set of cardinality $2^k$.
\end{dfn}
For any quadratic nonresidue $z$, the relation $z^{2^{k-1}n}\equiv -1\bmod p$ is true. This implies that for any $0\ne m<2^k$, the relation $z^{mn}\equiv -1\bmod p$ is true. Thus, set $\mu_{2^k}$, which is the 2-Sylow subgroup of the group of units $\F_p^{\times}$ in a finite field, is independent of the generator $z$, and contains precisely $2^k$ elements. 

\begin{dfn} \label{dfn8181.200E} \hypertarget{dfn8181.200E}  Let $p = 2^kn + 1$ be prime, where $k\geq3$ and $n\geq1$ is odd, and let $a$ be a quadratic residue modulo $p$. The root sequence mod $p$ is defined by
	$$\mathscr{R}_p=\{a^n,a^{2n},a^{3n},\ldots,a^{(2^{k}-1)n}\}.$$
\end{dfn}
The group of unity $\mu_{2^k}$  and the root sequence $\mathscr{R}_p$ generated by an arbitrary quadratic residue $z\in\F_p$ provide a complete set of information about the form of the square root of $a\in \F_p$. The longest possible root subsequence of even powers

\begin{equation}\label{eq8181.200d}
	a^n,\quad a^{2n},\quad a^{2^2n},\quad\ldots,\quad a^{2^{k-3}n}=\pm i,\quad a^{2^{k-2}n}=-1\quad a^{2^{k-1}n}=1.
\end{equation}

is related to the pseudo prime test to base $a\ne1$, consult the literature on pseudo primes. 
\begin{lem}\label{lem8181.100A} \hypertarget{lem8181.100} Let a and $z$ be quadratic residue and quadratic nonresidue mod $p = 2^kn + 1$, $n$ odd, respectively. Then the sequence of integers

\begin{equation}\label{eq8181.200f}
\omega_0=a^n,\quad \omega_1=\omega_0z^{n2^{k-m_0}},\quad  \omega_2=\omega_1z^{n2^{k-m_1}},\quad  \ldots, \quad \omega_{i+1}=\omega_iz^{n2^{k-m_{i}}}=1 ,
\end{equation}
where $m_i<k$ and the multiplicative order $\ord_p\omega_i=2^{m_i}$ decreases to $1$.
\end{lem} 
\begin{proof}[\textbf{Proof}] Since $a$ is a quadratic residue, and the element $\omega_0=a^n$ is a $2^k$ root of unity, the multiplicative order $\ord_p\omega_0=\ord_p a^n\mid 2^{k-1}$. Assume the maximal multiplicative order
\begin{equation}
	\omega_0^{ 2^{k-2}}=-1\quad \text{ and }\quad 	\ord_p\omega_0= 2^{m_0}= 2^{k-1}.
	\end{equation}
Since $0\ne z$ is a quadratic nonresidue and $z^{n2^{k-1}}=-1$ in $\F_p$, taking power yields
\begin{align}
\omega_1^{2^{m_0-1}}&=\left( \omega_0z^{n2^{k-m_0}}\right) ^{2^{m_0-1}}\\
&=\omega_0 ^{2^{m_0-1}}z^{n2^{k-1}}\nonumber\\
&=(-1)(-1)=1\nonumber.
\end{align}	
This implies that the multiplicative order $\ord_p\omega_1\mid 2^{m_0-1}= 2^{k-1-1}= 2^{k-2}$ decreases by at least one unit. Assume the maximal multiplicative order
\begin{equation}
	\omega_1^{ 2^{k-3}}=-1\quad \text{ and }\quad 	\ord_p\omega_1= 2^{m_0}= 2^{k-2}.
\end{equation}
Since $0\ne z$ is a quadratic nonresidue and $z^{n2^{k-1}}=-1$ in $\F_p$, taking power yields
\begin{align}
	\omega_2^{2^{m_1-1}}&=\left( \omega_1z^{n2^{k-m_1}}\right) ^{2^{m_1-1}}\\
	&=\omega_i ^{2^{m_1-1}}z^{n2^{k-1}}\nonumber\\
	&=(-1)(-1)=1\nonumber.
\end{align}	
This implies that the multiplicative order $\ord_p\omega_2\mid 2^{m_1-1}= 2^{k-2-1}= 2^{k-3}$ decreases by at least one unit. \\

Iterating this process by at most $i+1=k-1$ times. Assume the maximal multiplicative order
\begin{equation}
	\omega_i^{ 2^{1}}=-1\quad \text{ and }\quad 	\ord_p\omega_i= 2^{m_i}= 2^{2}.
\end{equation}
Since $0\ne z$ is a quadratic nonresidue and $z^{n2^{k-1}}=-1$ in $\F_p$, taking power yields
\begin{align}
	\omega_{i+1}^{2^{m_i-1}}&=\left( \omega_iz^{n2^{k-m_i}}\right) ^{2^{m_i-1}}\\
	&=\omega_0 ^{2^{m_i-1}}z^{n2^{k-1}}\nonumber\\
	&=(-1)(-1)=1\nonumber.
\end{align}	
This implies that the multiplicative order $\ord_p\omega_{i+1}\mid 2^{m_i-1}= 2^{k-(k-1)-1}= 1$, (it has reached a cycle of period 1). 
\end{proof}  
This routine is used in the calculation of the square root modulo a prime $p$. The algorithm, called the Tonelli-Shank algorithm, is widely available in the literature, see \cite{TA1891}, \cite[Section 3.5]{MV1997}. Similar proofs to the one above, and background details appear in paper in the literature, see \cite[Section 9.2]{KR1998}.	The general Tonelli-Shank algorithm is a random algorithm, see the discussion in \cite[p.\;99]{CP2005}. A result introduced in \hyperlink{thm1234.100}{Theorem} \ref{thm1234.100} computes a quadratic nonresidue modulo a prime $p$ in polynomial time. This innovation seems to change it into a deterministic algorithm.

\section{Square Roots Mod $p=8n+1$}\label{S8182}
A root of the equation $x^2-a\equiv  0 \bmod p$ such that $a^{(p-1)/2}\equiv 1 \bmod p$, is determined by a series of approximations. More precisely, the square root of any square $a\in \F_p$ is of the form $\sqrt{a}=\pm\mu a^{(n+1)/2}$ , where $\mu \in \mu_{2^k}$ is a unique root of unity.The number of
iterations in the algorithm used is mostly a function of the 2-adic valuation $v_2(p-1) = k$. Extensive details on these algorithms are available in references \cite{BH1999} to \cite{TS1999} and other sources. More general techniques  for computing roots of polynomials are introduced in \cite[Section 4.3]{LN1997}  and similar references.
\begin{thm}\label{thm8181.200} \hypertarget{thm8181.200} Let $p=2^kn+1$ be a prime and let $z\in\F_p$ be a quadratic nonresidue, where $k\geq3$ and $n\geq1$ is odd. Then, the square root of the a quadratic residue $0\ne a\in\F_p$ is given by 
$$\sqrt{a}=\pm a^{\frac{n+1}{2}}z^{n(2^{k-m_0-1}+2^{k-m_1-1}+2^{k-m_2-1}+\cdots+2^{k-m_i-1})},$$
where $m_j<k$.
\end{thm}
\begin{proof}[\textbf{Proof}] Let $\omega_0=a^n$ be a $2^k$ root of unity. Then the sequence of approximations by square elements
	\begin{align}\label{eq8181.200r}
		a\omega_0&=a^{n+1}=\left( a^{\frac{n+1}{2}}\right) ^2,
	\end{align}
	\begin{align}\label{eq8181.200s}
		a\omega_1&=a\omega_0z^{n2^{k-m_0}}\\&=a^{n+1}z^{n2^{k-m_0}}\nonumber\\
		&=\left( a^{\frac{n+1}{2}}z^{n2^{n(k-m_0-1)}}\right) ^2\nonumber,
	\end{align}	
	
	\begin{align}\label{eq8181.200t}
		a\omega_2&=a\omega_1z^{2^{k-m_1}}\\
		&=\left( a^{n+1}z^{2^{k-m_0}}\right) z^{2^{k-m_1}}\nonumber\\&=\left( a^{\frac{n+1}{2}}z^{n(2^{k-m_0-1}+2^{k-m_1-1})}  \right) ^2\nonumber,
	\end{align}	
	
	\begin{align}\label{eq8181.200u}
		a\omega_3&=a\omega_2z^{2^{k-m_2}}\\
		&=\left( a^{n+1}z^{n(2^{k-m_0-1}+2^{k-m_1-1})}\right) z^{2^{k-m_2}}\nonumber\\
		&=\left( a^{\frac{n+1}{2}}z^{n(2^{k-m_0-1}+2^{k-m_1-1})}  \right) ^2\nonumber,
	\end{align}	
	
	\hskip 1.5 in \ldots \hskip 1 in \ldots \hskip 1 in \ldots 
	
	\begin{align}\label{eq8181.200v}
		a\omega_{i+1}&=a\omega_iz^{2^{k-m_i}}\\
		&=\left(a^{n+1}z^{n(2^{k-m_0}+2^{k-m_1}+\cdots+2^{k-m_{i-1}})}\right) z^{2^{k-m_i}}\nonumber\\
		&=\left( a^{\frac{n+1}{2}}z^{n(2^{k-m_0-1}+2^{k-m_1-1}+2^{k-m_2-1}+\cdots+2^{k-m_i-1})}  \right) ^2\nonumber,
	\end{align}	
	
	converges to $a\in\F_p$.
\end{proof}

\begin{thm}\label{thm8181.200I} \hypertarget{thm8181.200I} Let $p=2^kn+1$ be a prime and let $z\in\F_p$ be a quadratic nonresidue, where $k\geq3$ and $n\geq1$ is odd. Then, there exists a unique integer $m\leq 2^{k-1}$, such that the square root of the quadratic residue $0\ne a\in\F_p$ is given by 
	
\begin{equation}\label{eq8181.200j}
\sqrt{a}=
\begin{cases}
\pm a^{\frac{n+1}{2}}z^{(2m+1)n}&\text{ if } \ord_p a^n=2^{k-1},\\[.3cm]
\pm a^{\frac{n+1}{2}}z^{2mn}&\text{ if } \ord_p a^n<2^{k-1}\nonumber.
\end{cases}
\end{equation}

\end{thm}
\begin{proof}[\textbf{Proof}] After about $k$ iterations, the sequence of approximations in \hyperlink{thm8181.200I}{Theorem} \ref{thm8181.200I} yields
	\begin{equation}\label{eq8181.200l}
a\omega_{i+1}=\omega_ia^{n+1}z^{n2^{k-m_i}}=\left( a^{\frac{n+1}{2}}z^{n(2^{k-m_0-1}+2^{k-m_1-1}+2^{k-m_2-1}+\cdots+2^{k-m_i-1})} \right) ^2.
	\end{equation}
Moreover the identity $\omega_{i+1}=\omega_iz^{n2^{k-m_i}}=1$ follows from   \hyperlink{lem8181.100}{Lemma} \ref{lem8181.100A}. This implies that
	\begin{align}\label{eq8181.200m}
	a&=\left( a^{\frac{n+1}{2}}z^{n(2^{k-m_0-1}+2^{k-m_1-1}+2^{k-m_2-1}+\cdots+2^{k-m_i-1})} \right) ^2\\
	&=\left( a^{\frac{n+1}{2}}z^{nr} \right) ^2\nonumber,
\end{align}
where the exponent
\begin{equation}\label{eq8181.200k}
	r=2^{k-m_0-1}+2^{k-m_1-1}+2^{k-m_2-1}+\cdots+2^{k-m_i-1}\leq 2^{k-1}
\end{equation} 
since $m_i\in\{0,1\}$. The parity of the integer $r\geq1$ depends on the multiplicative order $\ord_p a^n$. 
\end{proof}

Since the shortest sequence $ a^n = 1$ occurs if and only if the multiplicative order $\ord_p(a) = \text{odd}$, the probability that a quadratic residue has
odd order is $1 / 2^{k-1}$. Thus for $2^k\geq 8$, less than 1\% of the squares $a\in\F_p$ have square roots of the form $\sqrt{a}=\pm a^{\frac{n+1}{2}}z^{2mn}$ .

\begin{cor}\label{cor8181.200Q} \hypertarget{cor8181.200Q}
As $n \to\infty$, $k\geq3$ fixed, almost (in the sense of natural density) every quadratic residue has a square root of the
form $\sqrt{a}=\pm a^{\frac{n+1}{2}}z^{2mn}$, where $m \leq  2^{k-1}$.
\end{cor}

\begin{proof}[\textbf{Proof}] The exponent $e$ of the multiplier $z^{en}$ is odd if and only if the quadratic residue $a\in\F_p$ has multiplicative order $\ord_p(a) = 2^{k-1}$. But the
probability that $\ord_p(a) = 2^{k-1}$ is $2^{k-2}/n2^{k-1} = 1/2n$.
\end{proof}
Some statistical results are discussed in \cite{TS1999}.
\section{Square Roots Mod $p\ne8n+1$}\label{S8185}
The difficult square root calculations occurs only for primes in the arithmetic progression $ p\equiv 1 \bmod 8$. The square root calculations for the other 3 arithmetic progressions $p\equiv 3 \bmod 8$, $p\equiv 5 \bmod 8$ and $ p\equiv 7 \bmod 8$ are well known and are much simpler. These are included here for completeness.

\begin{thm}\label{thm8185.300} \hypertarget{thm8185.300} If $p=2^2n+3$ is a prime, then the square root of the a quadratic residue $0\ne a\in\F_p$ is given by 
	\begin{equation}\label{eq8185.300j}
		\sqrt{a}=
				\pm a^{n+1}.
	\end{equation}	
\end{thm}
\begin{proof}[\textbf{Proof}] Observe that $n+1=(p+1)/2=(p-1)/2+1$. Next, squaring both sides yields
\begin{equation}\label{eq8185.300l}
\left( 	\sqrt{a}\right) ^2=\left( 
\pm a^{n+1}\right) ^2=\left( 
	\pm a^{\frac{p+1}{4}}\right) ^2=a^{\frac{p-1}{2}}\cdot a=a,
\end{equation}		
since $a\ne0$ is a quadratic residue modulo $p=2^2n+3$.
\end{proof}

\begin{thm}\label{thm8185.350} \hypertarget{thm8185.350} If $p=2^3n+5$ is a prime, then the square root of the a quadratic residue $0\ne a\in\F_p$ is given by 
	\begin{equation}\label{eq8185.350k}
		\sqrt{a}=
		\begin{cases}
			\pm a^{n+1},&\text{ if }  a^{(p-1)/4}\equiv 1\bmod p,\\[.3cm]
			\pm 2^{2n+1}a^{n+1},&\text{ if }  a^{(p-1)/4}\equiv -1\bmod p\nonumber.
		\end{cases}
	\end{equation}	
\end{thm}
\begin{proof}[\textbf{Proof}] The first case for $a^{(p-1)/4}\equiv 1\bmod p$ is similar to \eqref{eq8185.300l}. The second case for $a^{(p-1)/4}\equiv -1\bmod p$ is slightly different. Take the quadratic nonresidue $z=2$ modulo $p=2^3n+5$, and $n=(p-5)/8$. Then
\begin{equation}\label{eq8185.300m}
	\left( 	\sqrt{a}\right) ^2=\left( \pm a^{n+1}2^{2n+1}\right) ^2=\left( 
	 a^{\frac{p-1}{4}}\cdot a\right) 2^{\frac{p-1}{2}}=(-1)a\cdot (-1)=a,
\end{equation}

since $2^{(p-1)/2}\equiv -1 \bmod p$.
\end{proof}	
In the second case, the form $	\sqrt{a}=\pm 2a(4a)^{(p-5)/8}$ is more efficient in numerical calculations. The earlier work on this case appears in \cite[p.\;219]{DL1966}.

\section{Polynomial Representations of the Square Roots}\label{S8183}
Any one-to-one function $f(x)$ on a finite field $\F_q$ has a polynomial representation $f(x) = a_dx^d + a_{d-1}x^{d-1} + \cdots + a_1x + a_0\in \F_p[x]$ of degree
$\deg(f) = d \leq p -2$, see \cite[Theorem 7.6]{LN1997}. As stated before the polynomial representation $f(x)\in \F_p[x]$ of the square root mod $p$ satisfies the
equation $\sqrt{x}\equiv \pm f (x)\bmod p$, where $x$ a quadratic residue. The principal root 
$\sqrt{x}\equiv f (x)\bmod p$ is a one-to-one function.
\begin{thm}\label{thm8183.500} \hypertarget{thm8183.500} Let $p=2^kn+1$ is a prime with $k\geq3$ and $n\geq1$ odd. Then there is a polynomial representation $f(x)\in \F_p[x]$ of the square root function of degree $\deg(f) = 2^{k-1}n- (n-
	1)/2$, and $2^k-1$ terms. Moreover it has the form 
	
\begin{equation}\label{eq8183.500j}
f(x)=\pm 2^{-(k-1)}x^{\frac{n+1}{2}}\left(c_dx^{dn}+ c_{d-1}x^{(d-1)n}+ c_{d-2}x^{(d-2)n}+\cdots+c_2x^{2n}+c_1x^{n}\right), \nonumber
\end{equation}		
where $d=2^{k-1}-1$.
\end{thm}

The polynomial $f(x)=f_k(x)$ is computed by considering all the possible sequences $a^n, a^{2n} ,a^{3n}, a^{4n}, a^{5n}, \ldots$ generated by an
arbitrary quadratic residue a modulo p and the sequence given in \eqref{eq8181.200f}. Each combination of the sequence is then mapped to
a unique term
\begin{equation}\label{eq8183.500l}
z^{a_1n}\left( 1\pm x^{a_2n}z^{a_3n}\right) \left( 1\pm x^{a_4n}z^{a_5n}\right) \cdots \left( 1\pm x^{a_{k-1}n}z^{a_kn}\right) ,
\end{equation}		
where $0\leq a_i < 2^k-1$. For example, the first term $z^{7n}(1- x^{4n})(1- x^{2n}z^{4n})(1- x^nz^{6n})$ in the polynomial $f_4(x)$ for the primes $p =
2^4n + 1$ corresponds to the longest sequence $a^n, a^{2n}, a^{3n}, a^{4n}, a^{5n}, \ldots, -1, 1$, and the last term $(1 + x^{4n})(1 + x^{2n})(1 + x^n)$
corresponds to the shortest sequence $a^n = 1$.\\

The first few of these polynomials are given below. The first two formulae are well known, but the next two are new.
\begin{enumerate}
\item  A polynomial for the subset of primes $p=2^2n+3$, $n\geq1$. 
\begin{align}	\sqrt{x}=
	\pm f_1(x)=\pm x^{n+1}.
\end{align}	
Here the degree is $\deg f_1=n+1=(p+1)/4$, see \hyperlink{thm8185.300}{Theorem} \ref{thm8185.300}. \\	

\item A polynomial for the subset of primes $p=2^3n+5$, $n\geq1$. 
\begin{align}	\sqrt{x}=
	\pm f_2(x)=\pm\frac{1}{2} x^{n+1}\left( \left( x^{2n+1}+1\right) -2^{2n+1}\left( x^{2n+1}-1\right) \right). 
\end{align}	
Here the degree is $\deg f_2=3n+2=(3p+1)/8$, this follows from \hyperlink{thm8185.350}{Theorem} \ref{thm8185.350}. \\

\item A polynomial for the subset of primes $p=2^3n+1$, $n$ odd. 
\begin{align}
\sqrt{x}&=
\pm f_3(x)\\
&=\pm \frac{1}{4}x^{\frac{n+1}{2}} \left[ z^{3n}(1-x^{2n} )(1-x^nz^{2n} )+ z^{2n}(1+x^{2n} )(1-x^n )\right .\nonumber\\
&\left . \hskip .75 in + z^{n}(1-x^{2n} )(1-x^nz^{2n} )+ (1+x^{2n} )(1+x^n)\right].\nonumber
\end{align}
Here the degree is $\deg f_3=3n+(n+1)/2=(7p+1)/16$, this follows from \hyperlink{thm8183.500}{Theorem} \ref{thm8183.500}.

\item  A polynomial for the subset of primes $p=2^4n+1$, $n$ odd. 
\begin{align}
	\sqrt{x}&=
	\pm f_4(x)\\
	&=\pm \frac{1}{8}x^{\frac{n+1}{2}} \left[ z^{7n}(1- x^{4n})(1- x^{2n}z^{4n})(1- x^nz^{6n})\right .\nonumber\\
	&\hskip .75 in \left.+z^{6n}(1+ x^{4n})(1- x^{2n})(1- x^{n}z^{4n})\right .\nonumber\\
	&\hskip .75 in\left.+z^{5n}(1-x^{4n})(1-x^n z^{2n})(1+ x^{2n}z^{4n})\right .\nonumber\\
	&\hskip .75 in\left.+z^{4n}(1+x^{4n})(1+x^{2n})(1- x^{n})\right .\nonumber\\
	&\hskip .75 in\left.+z^{3n}(1+x^{4n})(1-x^{2n} z^{4n})(1+ x^{n}z^{6n})\right .\nonumber\\
	&\hskip .75 in\left.+z^{2n}(1+x^{4n})(1-x^{2n})(1+ x^{n}z^{4n})\right .\nonumber\\
	&\hskip .75 in\left.+z^{n}(1-x^{4n})(1+x^{2n})(1+ x^{n}z^{4n})\right .\nonumber\\
	&\hskip .75 in\left.+(1+x^{4n})(1+x^{2n})(1+ x^{n})\right].\nonumber
\end{align}
Here the degree is $\deg f_4=4n+(n+1)/2=(15p+1)/32$, this follows from \hyperlink{thm8183.500}{Theorem} \ref{thm8183.500}.
\end{enumerate}

\section{The Least Quadratic Nonresidue Mod $p$}\label{S2424}
The least quadratic nonresidue $n(p)\in[2,p-1]$ modulo a prime $p$ has an old unconditional upper bound of $n(p)\ll p^{e^{-1/2}+\varepsilon}$ and a conjectured value of $n(p)\ll (\log x)^{1+\varepsilon}$, a recent survey of the literature appears in \cite{MT2021}. A new result  proves that there exists a quadratic nonresidue $z(p)\ll (\log p)^2$, this is computable in polynomial time.  

\begin{thm}\label{thm1234.100} \hypertarget{thm1234.100} Let $p$ be a large prime and let $n(p)$ denotes the least quadratic nonresidue modulo $p$. Then,
	\begin{equation}\label{eq1234.100}
		n(p)\ll (\log p)^{1+\varepsilon},    
	\end{equation}
\end{thm}
where $\varepsilon>0$ is an arbitrary small number.
\begin{proof}[\textbf{Proof}] The complete proof appears in \cite{CN2021}.
\end{proof}


\end{document}